\documentclass[12pt,a4paper,twoside]{article}  
\usepackage{amsmath,amsthm,amssymb}

\swapnumbers 
\theoremstyle{plain}

\newtheorem{theorem}{Theorem}

\newtheorem{lemmae}[theorem]{Lemma}

\newtheorem{corollary}[theorem]{Corollary}

\newtheorem{proposition}[theorem]{Proposition}

\theoremstyle{definition}
\newtheorem{definition}[theorem]{Definition}

\theoremstyle{remark}

\def \cz {\Bbb C}  
\DeclareMathOperator{\logdet}{logdet}
\DeclareMathOperator{\Vol}{Vol}
\DeclareMathOperator{\Tor}{Tor}

\begin{document}  

\title{$L^2$-torsion of hyperbolic manifolds}
\author{Eckehard Hess \and Thomas Schick}
\date{to appear in Manuscripta Mathematica}
\maketitle
\begin{center}\scriptsize
\begin{tabular}[l]{l@{\hspace{6mm}}l}
Eckehard Hess & Thomas Schick\\
Fachbereich Mathematik & Fachbereich Mathematik\\
Universit\"at Mainz & Universit\"at M\"unster\\
& Einsteinstr.\ 62\\
55099 Mainz, Germany & 48149 M\"unster, Germany\\
hess@topologie.mathematik.uni-mainz.de& thomas.schick@math.uni-muenster.de
\end{tabular}
\end{center}\vspace{2mm}

The $L^2$-torsion is an invariant defined for compact $L^2$-acyclic
manifolds of determinant class, for example odd  dimensional hyperbolic
manifolds. It was introduced by John Lott \cite{Lott} and Varghese Mathai
\cite{matthai} and computed for hyperbolic manifolds in low dimensions. Our
definition of the $L^2$-torsion coincides with that of John Lott, which is
twice the logarithm of that of Varghese Mathai.\par
In this paper you will find a proof of the fact that  the
$L^2$-torsion of hyperbolic manifolds of arbitrary odd dimension does not vanish.
This was conjectured by John Lott in \cite[p.484, Proposition 16 infra]{Lott}.
Some concrete values are computed and an estimate of their growth with the
dimension is given. The values we compute for dimensions 5 and 7 differ from
those published in \cite[Proposition 16]{Lott}. The result has been
independently achieved by both authors and will be part of the dissertation
of Eckehard Hess at the university of Mainz. For an introduction into
$L^2$-theory see \cite{lueck3}.\par
We are indebted to Prof.~Dr.~L\"uck, M\"unster, for
permanent support and encouragement.

\begin{definition}\label{def}
Following \cite[p.482]{Lott} we define the analytic $L^2$-torsion
of an $L^2$-acyclic
Riemannian $(d=2n+1)$-dimensional manifold of determinant class by
$$\Tor_{(2)}(M)=2\sum_{j=0}^n
(-1)^{j+1}\logdet_{G}(\triangle_j|)$$
Here $G$ is the fundamental group of $M$, $\triangle_j|$ is the Laplacian
restricted to coclosed forms on the universal covering $\tilde{M}$ and the
logarithm of the $G$-determinant is computed from the local trace of the
heat kernel as follows
\begin{eqnarray*}
\logdet_{G}(\triangle_j|) &=& \int_F \Bigg{\{}
\frac{d}{ds}\Bigg{|}_{s=0}\Bigg{[}\frac{1}{\Gamma(s)}\int_0^1
t^{s-1}\mbox{tr}_{\cz}e^{-t\triangle_j|}(x,x) dt\Bigg{]}\\
&& + \int_1^\infty
t^{-1}\mbox{tr}_{\cz}e^{-t\triangle_j|}(x,x) dt\Bigg{\}}dx
\end{eqnarray*}
Here $F$ is a fundamental domain of $M$ in $\tilde{M}$, the first integral
exists for $s$ sufficiently large and one has to take the meromorphic
extension at $0$. $M$ being of determinant class ensures the second integral
to converge.
\end{definition}

\begin{theorem}\label{hyperbo}
There is a constant $\alpha_d>0$, such that for every $(d=2n+1)$-dimensional
closed hyperbolic manifold
$$\Tor_{(2)}(M)=(-1)^n \alpha_d \Vol(M)$$
The values for $\alpha_d$ have been computed as follows.
Although exact values were computed for $d\leq 251$, we will give only exact
numbers for $d\leq11$ and approximate numbers for $d\leq39$:
$$\begin{array}{r|c|r}
d & \alpha_d & \approx \alpha_d\\
\hline
3 & \frac{1}{3\pi}&0,106103  \\[1mm]
5 & \frac{62}{45\pi^2}& 0,139598 \\[1mm]
7 & \frac{221}{35\pi^3}& 0,203645 \\[1mm]
9 & \frac{32204}{945\pi^4}& 0,349847 \\[1mm]
11 & \frac{1339661}{6237\pi^5}& 0,701891 \\[1mm]
\end{array}\hspace{8mm}
\begin{array}{r|r}
d & \approx\alpha_d\\
\hline
13 & 1,61885 \\
15 & 4,22925 \\
17 &12,3578 \\
19 & 39,9606 \\
21 & 141,729 \\
23 &  547,188\\
\end{array}\hspace{8mm}
\begin{array}{r|r}
25 &  2284,87\\
27 & 10261,5 \\
29 & 49326\\
33 & 252701 \\
35 & 1,37458 *10^6 \\
37 & 7,91236 *10^6 \\
39 &4,80523 *10^7\\
\end{array}$$

\end{theorem}

\begin{lemmae}
Let $\triangle_j|=\delta_jd_j|_{ker(\delta_{j+1})}$ be the Laplacian,
restricted to coclosed
$L^2$-forms on the $(d=2n+1)$-dimensional
hyperbolic space $H^d$. Then for every closed $d$-dimensional
hyperbolic manifold $M$ with fundamental group $G$ and
$j\ne n$
$$ \logdet_G(\triangle_j|)=\Vol(M) C
{ \left(
\begin{array}{c}
 2n \\
j
\end{array}
\right)}
\sum_{k=0}^{n}K^n_{k,j}(-1)^{k+1}\frac{\scriptstyle
2\pi}{\scriptstyle 2k+1}(n-j)^{2k+1}$$
$$\mbox{with
}C=\frac{(4\pi)^{-(n+\frac{1}{2})}}{\Gamma(n+\frac{1}{2})}\mbox{\hspace
{8mm} and \hspace{8mm}}a=n-j$$
Here $K^n_{k,j}$ is the coefficient of $\nu^{2k}$ in the polynomial
$$P(\nu):=\frac{\prod_{i=0}^{n}(\nu^2+i^2)}{\nu^2+(n-j)^2}\hspace{8mm}(*)
\label{kdkj}$$
Note that $P(\nu)$ indeed is a polynomial rather than a rational function, as
$|a|\in\{1...n\}$. In addition
$$ \logdet_G(\triangle_{n}|)=0$$
\end{lemmae}

\begin{proof}
Following \cite[prop.~15]{Lott} the local trace of the heat kernel of
$\triangle_j|$ is
$$\mbox{tr}_{\cz} e^{-t\triangle_j|}(x,x)=C
{ \left(
\begin{array}{c}
 2n \\
j
\end{array}
\right)}
\int_{-\infty}^{\infty}e^{-t(\nu^2+a^2)}\frac{\prod_{k=0}^{n}
(\nu^2+k^2)}{\nu^2+a^2}d\nu$$
According to the above remark  let
$$\frac{\prod_{k=0}^{n}(\nu^2+k^2)}{\nu^2+a^2}=\sum_{k=0}^{n}K^n_{k,j}\nu^{2k
}$$
Evaluation of the above integral yields
$$\mbox{tr}_{\cz} e^{-t\triangle_j}(x,x)= C
{ \left(
\begin{array}{c}
 2n \\
j
\end{array}
\right)}
\sum_{k=0}^{n}K^n_{k,j}e^{-ta^2}t^{-k-\frac{1}{2}}
\Gamma\Big{(}k+ \frac{\scriptstyle 1}{\scriptstyle 2}\Big{)}$$
Now we have to compute
\begin{eqnarray*}
L_j & := &\frac{\logdet(\triangle_j|)}{\Vol(M)}\\
& = &\frac{d}{ds}\Bigg{|}_{s=0}\Bigg{[}\frac{1}{\Gamma(s)}\int_0^1 t^{s-1}
\mbox{tr}_{\cz}e^{-t\triangle_j}(x,x)dt\Bigg{]} + \int_1^\infty
t^{-1}\mbox{tr}_{\cz}e^{-t\triangle_j}(x,x) dt
\end{eqnarray*}

John Lott showed in \cite[Lemma 13, p.481]{Lott}

$$L_{n} =  0$$
For $j\ne n$, that is $a=n-j\ne 0$

\begin{eqnarray*}
L_j & = & C { \left(
\begin{array}{c}
 2n \\
j
\end{array}
\right)}
\sum_{k=0}^{n} K^n_{k,j}
\Gamma\Big{(}k+\frac{\scriptstyle1}{\scriptstyle2}\Big{)}
\frac{d}{ds}\Bigg{|}_{s=0}\Bigg{(}\frac{1}{\Gamma(s)}\underbrace{
\int_0^\infty
e^{-ta^2}t^{s-k-\frac{3}{2}}dt}_{:=J}\Bigg{)}
\end{eqnarray*}

$J$ exists for $s$ sufficiently large. Its meromorphic extension leads to
\begin{eqnarray*}
L_j & = &  C { \left(
\begin{array}{c}
 2n \\
j
\end{array}
\right)}
\sum_{k=0}^{n} K^n_{k,j}
\Gamma\Big{(}k+\frac{\scriptstyle1}{\scriptstyle2}\Big{)}
\Gamma\Big{(}-k-\frac{\scriptstyle1}{\scriptstyle2}\Big{)} a^{2k+1}\\
& = & C { \left(
\begin{array}{c}
 2n \\
j
\end{array}
\right)}
\sum_{k=0}^{n} K^n_{k,j} (-1)^{k+1}\frac{\scriptstyle
2\pi}{\scriptstyle 2k+1}a^{2k+1}
\end{eqnarray*}
\end{proof}

\begin{corollary}
For any closed hyperbolic manifold of dimension $d=2n+1$ we have by
Definition \ref{def}
$$\frac{\Tor_{(2)}(M)}{\Vol(M)}=2 \sum_{j=0}^{n-1} (-1)^{j+1} C
{ \left(
\begin{array}{c}
 2n \\
j
\end{array}
\right)}
\sum_{k=0}^{n} K^n_{k,j}
(-1)^{k+1}\frac{\scriptstyle2\pi}{\scriptstyle 2k+1} (n-j)^{2k+1}$$
\end{corollary}
The numerical values were computed using this fomula and Mathematica.

\begin{lemmae}
Let $M$ be a closed $(d=2n+1)$-dimensional manifold. Then
$$(-1)^{j+1} \logdet_G(\triangle_j|)=(-1)^n|
\logdet_G(\triangle_j|)|$$
with
$$| \logdet_G(\triangle_j|)|>0 \mbox{ for } j\not=n$$
In particular
$$(-1)^n\Tor_{(2)}(M)>0$$
\end{lemmae}

\begin{proof}
Let $j\not=n$. Then we have
\begin{eqnarray*}
L_j  & = &- 2 \pi C { \left(
\begin{array}{c}
 2n \\
j
\end{array}
\right)}
\sum_{k=0}^{n} K^n_{k,j} (-1)^k
   \frac{\scriptstyle 1}{\scriptstyle 2k+1}a^{2k+1}\\
 & = &    - 2\pi C { \left(
\begin{array}{c}
 2n \\
j
\end{array}
\right)}
\int_0^a \sum_{k=0}^{n} K^n_{k,j} (ix)^{2k}dx\\
\end{eqnarray*}
Using the definition $(*)$ of the coefficients $K^n_{k,j}$ in Lemma
\ref{kdkj} one
gets
\begin{eqnarray*}
L_j & = &  - 2 \pi C { \left(
\begin{array}{c}
 2n \\
j
\end{array}
\right)}
\int_0^a \frac{\prod_{k=0}^{n}(k^2-x^2)}{a^2-x^2}dx\\
\end{eqnarray*}
One has
\begin{eqnarray*}
 & & \int_0^a \frac{\prod_{k=0}^{n}(k^2-x^2)}{a^2-x^2}dx\\
& = & (-1)^{n+1}\int_0^a
\frac{x}{(a+x)(a-x)}\prod_{k=-n}^{n}(x+k)dx\\
& = & (-1)^{n+1}\sum_{r=0}^{a-1}\int_0^1 f_r(t)dt
\end{eqnarray*}
where for   $t\in]0,1[\;\;,\;r \in \{0,..,a-1\}$ we define
\begin{eqnarray*}
f_r(t) & = &\underbrace{\frac{t+r}{(a+t+r)(a-t-r)}}_{>0}
\prod_{k=-n+r}^{n+r}\underbrace{(t+k)}_{\begin{array}{c}
<0 \text{
for }k<0 \\ >0 \text{ otherwise}\end{array} } \\
 & = & (-1)^{n-r}|f_r(t)|
\end{eqnarray*}

For $t \in\; ]0,1[$ and $0\leq r <r+1 \leq a-1$ one computes
$$\left| \frac{f_{r+1}(t)}{f_r(t)}\right| >1$$

Hence
$$\int_0^1|f_{r+1}(t)|dt\geq \int_0^1|f_r(t)|dt$$
Now the sum
$$L_j= - 2\pi C { \left(
\begin{array}{c}
 2n \\
j
\end{array}
\right)}
\sum_{r=0}^{a-1}(-1)^{r+1}\int_0^1 |f_r(t)|dt$$
is an alternating sum and the absolute values of the summands are stricly
increasing. So it is not $0$ and the sign is that of the last summand. One
concludes
$$ \logdet_G(\triangle_j|)=(-1)^{n-j-1}|
\logdet_G(\triangle_j|)| $$
This also finishes the proof of Theorem \ref{hyperbo}.
\end{proof}

\begin{proposition}
The constants $\alpha_d$ of Theorem \ref{hyperbo} strictly increase and
$$\alpha_{2n+1}\ge \frac{n}{2\pi}\alpha_{2n-1}$$
In particular
$$\alpha_{2n+1}\ge\frac{2}{3}\frac{n!}{(2\pi)^n}$$
\end{proposition}

\begin{proof}
An elementary computation shows
$$ \left|\frac{f_{a-1}(t)}{f_{a-2}(t)}\right| \geq 2$$
Now one has
\begin{eqnarray*}
\left|\int_0^a \frac{\prod_{k=0}^n(k^2-x^2)}{a^2-x^2}dx\right|
&\geq& \frac{1}{2}\int_{a-1}^a
\left|\frac{\prod_{k=0}^n(k^2-x^2)}{a^2-x^2}\right|dx
\end{eqnarray*}
and
\begin{eqnarray*}
\alpha_{2n+1}&\geq&
2\pi\sum_{j=0}^{n-1}
\frac{(4\pi)^{-(n+\frac{1}{2})}}{\Gamma\Big{(}n+\frac{1}{2}\Big{)}}
\left({
\begin{array}{c}
2n\\j
\end{array}}\right)
\int_{n-j-1}^{n-j}\left|\frac{\prod_{k=0}^{n-1}(k^2-x^2)}{(n-j)^2-x^2}\right
| \underbrace{(n^2-x^2)}_{\geq(2n-j)j}dx\\
&\geq&
2\pi\sum_{j=1}^{n-1}
\frac{(4\pi)^{-(n+\frac{1}{2})}}{\Gamma\Big{(}n+\frac{1}{2}\Big{)}}
2n(2n-1)\left(
{ \begin{array}{c}
2n-2\\j-1
\end{array}}\right)
\int_{n-j-1}^{n-j}\left|\frac{\prod_{k=0}^{n-1}(k^2-x^2)}{(n-j)^2-x^2}\right
| dx\\
&=&
4\pi\sum_{j=1}^{n-1} \frac{2n}{4\pi}
 \frac{(4\pi)^{-(n-1+\frac{1}{2})}}{\Gamma\Big{(}n-1+\frac{1}{2}\Big{)}}
\left(
{ \begin{array}{c}
2n-2\\j-1
\end{array}}\right)
\int_{n-j-1}^{n-j}\left|\frac{\prod_{k=0}^{n-1}(k^2-x^2)}{(n-j)^2-x^2}\right
| dx\\
 &\geq&
\frac{n}{2\pi}
\;4\pi\sum_{l=0}^{n-2}
 \frac{(4\pi)^{-(n-1+\frac{1}{2})}}{\Gamma\Big{(}n-1+\frac{1}{2}\Big{)}}
\left(
{ \begin{array}{c}
2(n-1)\\l
\end{array}}\right)
\int_{n-1-l-1}^{n-1-l}\left|\frac{\prod_{k=0}^{n-1}(k^2-x^2)}{(n-1-l)^2-x^2}
\right| dx\\
&\geq&
\frac{n}{2\pi}
\;4\pi\sum_{l=0}^{n-2}
 \frac{(4\pi)^{-(n-1+\frac{1}{2})}}{\Gamma\Big{(}n-1+\frac{1}{2}\Big{)}}
\left(
{ \begin{array}{c}
2(n-1)\\l
\end{array}}\right)
\left|\int_{0}^{n-1-l}\frac{\prod_{k=0}^{n-1}(k^2-x^2)}{(n-1-l)^2-x^2}
dx\right|\\
&=& \underbrace{\frac{n}{2\pi}}_{\geq 1 \mbox{ for }n\geq 7}
\alpha_{2n-1}
\end{eqnarray*}
For $n\leq 7$ the growth follows from the table.
\end{proof}

\end{document}